

\input amstex
\documentstyle {amsppt}
\document
\centerline{\bf GENERATING FUNCTIONS IN ALGEBRAIC GEOMETRY}

\smallskip

\centerline{\bf AND SUMS OVER TREES}

\medskip

\centerline{Yu.I.Manin}

\smallskip

\centerline{\it Max--Planck--Institut f\"ur Mathematik, Bonn, Germany}
\document
\TagsOnRight
\nologo
\magnification=1200
\NoBlackBoxes
\bigskip

\centerline{\bf \S 0. Summary of results}

\medskip

{\bf 0.1. Introduction.} In this paper we adress the following three
problems.

\smallskip

A. Calculate the Betti numbers and Euler characteristics of moduli spaces
$\overline{M}_{0,n}$ of stable $n$--pointed curves of genus zero
(see e.g. [Ke]), or rather an appropriate generating function
for these numbers.

\smallskip

B. The same for the space $X[n]$, a natural compactification
of the space of $n$ pairwise distinct labelled points on a
non--singular compact algebraic variety $X$ constructed for
$\roman{dim}\ X=1$ in [BG] and in general in [FMPh].
(Beilinson and Ginzburg called this space ``Resolution
of Diagonals'', Fulton and MacPherson use the term
``Configuration Spaces'').

\smallskip

C. Calculate the contribution of multiple coverings in the problem
of counting rational curves on Calabi--Yau threefolds
(see [AM], [Ko], and more detailed explanations below).

\smallskip

All these problems are united by the fact that available
algebro--geometric information allows us to represent the corresponding
numbers as a sum over trees with markings. M. Kontsevich
in [Ko] invoked a general formula of perturbation theory
in order to reduce the calculation of the relevant
generating functions to the problem of finding the critical value
of an appropriate formal potential. We solve problems A and B
by applying this formalism in a simpler geometric context
than that of [Ko]. Problem C is taken from [Ko]; we were
able to directly complete Kontsevich's calculation in this case
and obtain a simple closed answer.

\smallskip

We will now describe our results (0.3---0.5) and technique
(0.6) in some detail.

\medskip

{\bf 0.2. General setup.} Let $Y$ be an algebraic variety
over $\bold{C}$, possibly non--smooth and non--compact.
Following [FMPh] we denote by $P_Y(q)$ the virtual Poincar\'e
polynomial of $Y$ which is uniquely defined by the following
properties.

\smallskip

a). If $Y$ is smooth and compact, then
$$P_Y(q)=\sum_{j}\roman{dim}\ H^j(Y)q^j. \eqno{(0.1)}$$
In particular
$$\chi (Y)=P_Y(-1). \eqno{(0.2)}$$

\smallskip

b). If $X=\coprod_iX_i$ is a finite union of pairwise disjoint
locally closed strata, then
$$P_Y(q)=\sum_iP_{Y_i}(q). \eqno{(0.3)}$$

\smallskip

c). $P_{Y\times Z}(q)=P_Y(q)P_Z(q).$ It follows that
if $Y$ is a fibration over base $B$ with fiber $F$
locally trivial in Zariski topology, then $P_Y(q)=P_B(q)P_F(q).$

\smallskip

A definition of $P_Y(q)$ can be given using the weight filtration on the
cohomology with compact support:
$$P_Y(q)=\sum_{i,j}(-1)^{i+j}\roman{dim}\ (\roman{gr}_W^jH_c^i(Y,\bold{Q}))
q^j. \eqno{(0.4)}$$

\smallskip

We apply the additivity formula (0.3) to the strata of the
natural stratifications of $\overline{M}_{0,n}$ and $X[n]$
in Problems A, B. These strata can be indexed by marked trees
describing various coalescing patterns of $n$--point configurations.

\smallskip

In [Ko], the role of $Y$ is played by a compactification $M(W)$
of the space of parametrised rational curves in some manifold $W$.
The relevant trees describe Gromov type degenerations of these
curves. Kontsevich calculates certain Chern numbers of
vector bundles over $M(W)$ and uses Bott's fixed point formula
instead of (0.3) in order to represent them as a sum of
local contributions. To make Bott's formula applicable,
Kontsevich assumes that $W$ is endowed with a torus action
and lifts this action to $M(W).$ (Actually, his $M(W)$ is not
a manifold but a smooth stack).

\medskip

{\bf 0.3. Moduli spaces.} We put
$$\varphi (q,t):=t+\sum_{n=2}^{\infty}P_{\overline{M}_{0,n+1}}(q)
\frac{t^n}{n!}\in\bold{Q}[q][[t]], \eqno{(0.5)}$$
$$\chi (t):=\varphi (-1,t)=t+\sum_{n=2}^{\infty}
\chi (\overline{M}_{0,n+1})\frac{t^n}{n!}\in\bold{Q}[[t]]. \eqno{(0.6)}$$

\smallskip

\proclaim{\quad 0.3.1. Theorem} a). $\varphi (q,t)$ is the unique
root in $t+t^2\bold{Q}[q][[t]]$ of any one of the following
functional/differential equations in $t$ with parameter $q:$
$$(1+\varphi )^{q^2}=q^4\varphi -q^2(q^2-1)t+1, \eqno{(0.7)}$$
$$(1+q^2t-q^2\varphi )\varphi_t=1+\varphi . \eqno{(0.8)}$$

b). $\chi$ is the unique root in $t+t^2\bold{Q}[[t]]$ of any one
of the similar equations
$$(1+\chi )\ \roman{log}(1+\chi )=2\chi -t, \eqno{(0.9)}$$
$$(1+t-\chi )\chi_t=1+\chi . \eqno{(0.10)}$$
\endproclaim

Equations (0.8) and (0.10) are equivalent to the following recursive formulas
for the Poincar\'e
polynomials. Put $p_n=p_n(q)=P_{\overline{M}_{0,n+1}}/n!$.

\smallskip

\proclaim{\quad 0.3.2. Corollary} We have for $n\ge 1$:
$$
(n+1)p_{n+1}=p_n+q^2\sum_{{i+j=n+1}\atop{i\ge 2}}jp_ip_j, \eqno{(0.11)}
$$
$$
P_{\overline{M}_{0,n+2}}(q)=P_{\overline{M}_{0,n+1}}(q)+
q^2\sum_{{i+j=n+1}\atop{i\ge 2}}{{n}\choose{i}}P_{\overline{M}_{0,i+1}}(q)
P_{\overline{M}_{0,j+1}}(q). \eqno{(0.12)}
$$
\endproclaim

\smallskip

One can compare (0.11) with recursive formulas in [Ke], p. 550.

\smallskip

{}From (0.10) one sees that the function inverse to $\chi$ has a critical
point at $t=e-2$. Don Zagier has shown me how to derive from this
the following asymptotical formula:
$$
\chi (\overline{M}_{0,n+1})\cong\frac{1}{\sqrt{n}}\left(\frac{n}{e^2-2e}
\right)^{n-\frac{1}{2}}.
$$

We will prove Theorem 0.3.1 in \S 1. We will also discuss the ramification
properties of $\varphi$ as a function of $t$ for $q^2\ne 1.$

\medskip

{\bf 0.4. Configuration spaces.} For a compact smooth algebraic
manifold $X$ of dimension $m$, set
$$
\psi_X(q,t)=1+\sum_{n\ge 1}P_{X[n]}(q)\frac{t^n}{n!}
\in\bold{Q}[q][[t]], \eqno{(0.13)}
$$
$$
\chi_X(t)=\psi_X(-1,t)=1+\sum_{n\ge 1}\chi (X[n])\frac{t^n}{n!}
\in\bold{Q}[[t]]. \eqno{(0.14)}
$$
Put also
$$
\kappa_m=\frac{q^{2m}-1}{q^2-1}=P_{\bold{P}^{m-1}}(q).
$$

\smallskip

\proclaim{\quad 0.4.1. Theorem} Denote by $y^0=y^0(q,t)$
the unique root in $t+t^2\bold{Q}[q^2][[t]]$ of any one
of the following equations:
$$
\kappa_m(1+y^0)^{q^{2m}}=q^{2m}(q^{2m}+\kappa_m-1)y^0
-q^{2m}(q^{2m}-1)t+\kappa_m, \eqno{(0.15)}
$$
$$
\left[ q^{2m}t+1-(q^{2m}-1+\kappa_m)y^0\right]y^0_t=1+y^0.
\eqno{(0.16)}
$$
Then we have in $\bold{Q}[q][[t]]$:
$$
\psi_X(q,t)=(1+y^0)^{P_X(q)}.
\eqno{(0.17)}
$$
\endproclaim

\smallskip

\proclaim{\quad 0.4.2. Theorem} Denote by $\eta =\eta (t)$
the unique root in $t+t^2\bold{Q}[[t]]$ of any one of the
following equations:
$$
m(1+\eta )\ \roman{log}(1+\eta )=(m+1)\eta -t, \eqno{(0.18)}
$$
$$
(t+1-m\eta )\eta_t =1+\eta . \eqno{(0.19)}
$$
Then we have in $\bold{Q}[[t]]$:
$$
\chi_X(t)=(1+\eta)^{\chi (X)}.
\eqno{(0.20)}
$$
\endproclaim

\smallskip

Theorems 0.4.1 and 0.4.2 are proved in \S 2.

\smallskip

I am grateful to C. Soul\'e who remarked that (0.17) follows from a
less neat identity which I deduced initially.
He has also informed me that he and H. Gillet constructed a map
$X\mapsto [h^*(X)]$ from varieties to the $K_0$--ring of Grothendieck's
motives having all the formal properties of the virtual Poincar\'e polynomial.
We can more or less mechanically use it in all our constructions;
in particular, $q^2$ will be replaced by Tate's motive $[h^2(\bold{P}^1)].$

\smallskip

For the reader's convenience, we list the first terms of the generating
series we have considered:
$$
\varphi (q,t)=t+\frac{t^2}{2!}+\frac{t^3}{3!}(q^2+1)+
\frac{t^4}{4!}(q^4+5q^2+1)+\frac{t^5}{5!}(q^6+16q^4+16q^2+1)+
$$
$$
\frac{t^6}{6!}(q^8+42q^6+127q^4+42q^2+1)+\dots ,
$$
$$
P^{-1}\varphi_X(q,t)=t(\kappa_m+P-1)+
$$
$$
\frac{t^2}{2!}
\left[(P-1)(P-2)+\kappa_m(q^{2m}-2)+3(P-1)\kappa_m+3\kappa_m^2\right]+
$$
$$
\frac{t^3}{3!}
[P^3-6P^2+11P-6+\kappa_m(6P^2-26P+26+4Pq^{2m}-9q^{2m}+q^{4m})+
$$
$$
\kappa_m^2(15P+10q^{2m}-35)+15\kappa_m^3] +\dots
$$
where we put $P=P_X(q).$

\medskip

{\bf 0.5. Multiple coverings.} Consider the following general problem
of enumerative geometry.

\smallskip

{\bf Problem $P_{g,k}(X,\beta ,\Cal{I})$.} {\it Given a projective algebraic
manifold
$X$, find the number of parametrised algebraic curves of genus $g$ in $X$,
in the homology class $\beta$,
with $k$ marked points, satisfying some incidence conditions $\Cal{I}.$}

\smallskip

Notice that in this vaguely stated problem we implicitly assume that
the number of solutions is only ``virtually'' finite, and look for
the number of virtual solutions.

\smallskip

In [Ko], Maxim Kontsevich suggested a general scheme allowing him
to simultaneously define this number for a wide class of problems
and to calculate it in many cases using Bott's residue formula.
In the three examples he considered in full detail we have $X=\bold{P}^n$
for some $n,\ g=0,$ and $\beta$ is $d[\bold{P}^1]$ for some $d\ge 1.$
The remaining data is as follows.

\smallskip

(i)  $n=2:\ X=\bold{P}^2,\ k=3d-1.$
The problem is to find the number of rational curves of degree $d$
in $\bold{P}^2$ passing through $3d-1$ points in general positions.

\smallskip

(ii)  $n=4:\ X=\bold{P}^4,\ k=0.$ The problem is to find the number
of rational curves of degree $d$ lying in a quintic hypersurface $V.$

\smallskip

(iii)  $n=1:\ X=\bold{P}^1,\ k=0.$ Here we additionally assume that $X$
is a rational curve embedded in the quintic threefold
(or a more general Calabi--Yau threefold) with normal sheaf
$\Cal{O}(-1)\oplus\Cal{O}(-1)$, and the problem is to calculate
the contribution of maps of degree $d,$ $\bold{P}^1\to X,$
to the number of solutions of problem (ii).

\smallskip

Using a different definition of the last contribution which we denote $m_d$
P. Aspinwall and D. Morrison [AM] calculated it and confirmed an
earlier prediction by P. Candelas et al.

\smallskip

In this note we show that Kontsevich's formula gives the same
answer:

\smallskip

\proclaim{\quad 0.5.1. Theorem} $m_d=d^{-3}.$
\endproclaim

\medskip

{\bf 0.6. Summation over trees.} A {\it tree} $\tau$ here is a finite
connected simply connected CW--complex. We denote by $V_\tau$
the set of its vertices, $E_\tau$ the set of its edges. Valency
$|v|$ of a vertex $v\in V_{\tau}$ is the number of edges
adjoining $v$. A {\it flag} of $\tau$ is a pair $(v,e)$ where
$v$ is a vertex, and $e$ is an adjoining edge.

\smallskip

A {\it marking} of a tree $\tau$ is a vaguely defined notion.
It may consist of a family of marks of given type(s) put onto
vertices, edges, flags, and satisfying certain restrictions.
Below we will describe pecisely a family of markings which we will call
{\it standard} ones.

\smallskip

The generating functions $\varphi$ studied above and in [Ko] are calculated
in three steps.

\smallskip

{\it STEP 1.} Represent $\varphi$ as an (infinite) sum of
certain weights $w_{\varphi}(\tau ,\mu )$ taken over isomorphism
classes of marked trees $(\tau ,\mu )$:
$$
\varphi =\sum_{(\tau ,\mu )/(iso)}w_{\varphi }(\tau ,\mu ). \eqno{(0.21)}
$$
This stage involves a combinatorial encoding of the raw
algebro--geometric data, determining type of marking and weights.

\smallskip

{\it STEP 2.} Try to rewrite (0.21) in a standard form of the following
type. Choose a set $A$ (finite or countable) and a family of
{\it symmetric} tensors indexed by $A$: $g^{ab},\ a,b\in A;$
$C_{a_1,\dots ,a_k},\ a_i\in A,\ k\ge 1.$ The coordinates
$g^{ab},\ C_{a_1,\dots ,a_k}$ must be elements of a topological
commutative ring.

\smallskip

The {\it standard marking} corresponding to this data
is a map $f:\ F_{\tau}\to A.$

\smallskip

The {\it standard weight} of a marked tree $(\tau ,f)$
corresponding to this data is
$$
w(\tau ,f):=\frac{1}{|\roman{Aut}\ \tau|}\prod_{\alpha\in E_{\tau}}
f^{(\partial \alpha)}\prod_{v\in V_{\tau}}C_{f(\sigma v)}.
\eqno{(0.22)}
$$
Here we use the following notation. For an edge $\alpha ,$
$\partial \alpha$ denotes the set of two flags of this edge,
and $f(\partial \alpha )$ is the set of two marks $(a,b)$
put on these flags by $f.$ Similarly, for a vertex $v$,
$\sigma v$ denotes the set of all flags containing $v$, and
$f(\sigma v)$ is the respective family of marks.

\smallskip

Finally, the standard sum over trees, or in physics speak,
a {\it partition function} is
$$
Z:=\sum_{\tau /(iso)}\sum_{f:F_{\tau}\to A}w(\tau ,f).
\eqno{(0.23)}
$$

The passage from (0.21) to (0.23) is not completely automatic
and indeed not always possible. Luckily, in can be made
for all the problems discussed in [Ko] and here. I cannot
explain conceptually why this is so. In particular, the factor
$1/|\roman{Aut}\ \tau |$ in the Problems A, B, resp. C,
occurs for different geometric reasons.

\smallskip

If we managed to represent (0.21) in the form (0.23), then we can try
to complete the calculation of $\varphi =Z$ with the help
of the following identity.

\smallskip

Assume that the matrix $(g^{ab})$ has an inverse matrix $(g_{ab}).$

\smallskip

{\it STEP 3.} Consider an auxiliary family of independent variables
(fields) $\varphi =\{\varphi_a\ |\ a\in A\} .$ Construct
the formal function (potential)
$$
S(\varphi )=-\sum_{a,b\in A} g_{ab}\frac{\varphi_a \varphi_b}{2}+
\sum_{k\ge 1,a_i\in A}\frac{1}{k!}C_{a_1,\dots a_k}
\varphi_{a_1}\dots \varphi_{a_k}.
\eqno{(0.24)}
$$
Denote by $\varphi^0=\{\varphi^0_a\ |\ a\in A\}$ an appropriate
critical point of $S(\varphi )$ that is, a solution of equations
$\frac{\partial S}{\partial \varphi_a}|_{\varphi =\varphi^0}=0,\
a\in A.$

\smallskip

\proclaim{\quad 0.6.1. Claim}
$$
Z=S^{crit}=S(\varphi^0). \eqno{(0.25)}
$$
\endproclaim

\smallskip

This remains a ``physical'' statement until we specify the relevant
topological ring containing $g$ and $C$, prove the existence and
uniqueness of $\varphi^0$, and the convergence of $S(\varphi^0).$
(See [Ko] for the standard physical argument ``proving'' 0.6.1).
For example, considering $(g^{ab},\ C_{a_1,\dots ,a_k})$ as
independent formal variables, one can treat (0.22) as a formal
series in these variables, and prove (0.6.1) as an identity
in a localization of this ring.

\smallskip

Anyway, STEP 3 involves three calculational difficulties.

\smallskip

a). We must be able to sum $S(\varphi ).$ In our problems A,B
this is easy. In [Ko],  a partial success is achieved,
reducing $S(\varphi )$ to a new potential which is quadratic
in $\varphi_a$ but highly non--linear in a finite set of new
auxiliary variables.

\smallskip

b). We must be able to solve $dS=0$ and to find $\varphi^0.$

\smallskip

c). We must be able to calculate $S(\varphi^0).$

\smallskip

The following trick, also well known to physicists, will allow us
in certain cases to avoid the last unpleasant calculation.

\smallskip

We will deform the data $(g^{ab},\ C_{a_1,\dots ,a_k})$ by
introducing independent parameters $t=\{ t_a\ |\ a\in A\}$
and replacing $C_a$ by $t_aC_a$. The rest of the data $A,\
g^{ab},\ C_{a_1,\dots ,a_k}$ for $k\ge 2$ remains unchanged.
Let $Z^t, S^t,\varphi^{0t}$ be respectively
the deformed partition function, potential, and the critical point.
Then we have
\proclaim{\quad 0.6.2. Claim} For all $a\in A,$ we have
$$
\frac{\partial Z^t}{\partial t_a}=C_a\varphi_a^{0t}.
\eqno{(0.26)}
$$
\endproclaim

\smallskip

{}From the view point of generating functions, we lose no information
replacing (0.25) by (0.26).

\smallskip

To deduce (0.26) from (0.25), one applies Claim 0.6.1 to $Z^t$
and differentiates in $t$:
$$
\frac{\partial Z^t}{\partial t_a}=
\frac{\partial}{\partial t_a}\left( S^t(\varphi )|_{\varphi =\varphi^{0t}})
\right) =
$$
$$
\sum_b\frac{\partial S^t(\varphi )}{\partial \varphi_b}
|_{\varphi =\varphi^{0t}}\frac{\partial \varphi_b^{0t}}{\partial t_a}
+\frac{\partial S^t}{\partial t_a}|_{\varphi =\varphi^{0t}_a}=
C_a\varphi_a^{0t}
$$
because $S^t$ depends on $t$ only via linear terms $\sum t_aC_a\varphi_a.$

\smallskip

On the other hand, to prove (0.26) in a formal context, one can
totally bypass Claim 0.6.1 and simply apply a universal inversion
formula to the formal map $(\varphi_a)\mapsto (\partial S^t/\partial
\varphi_a)$ giving simultaneously existence, uniqueness, and
expression for $\varphi^{0t}$ as a sum over trees. Such inversion
formulas are classical. The version closest to our needs is given
in [GK]; the only difference is that $\partial S^t/\partial \varphi_a$
at $0$ does not vanish. We leave details to the reader.

\smallskip

Functional equations (0.7), (0.9), (0.15), (0.18) are essentially
relations for coordinates of the critical point. Differential
equations are obtained from them by differentiating in $t$.

\smallskip

{\it Acknowledgements.} I am grateful to M. Kontsevich for many
enlightening explanations, and to Don Zagier for teaching me PARI.
 After this work was written,
I learned that E. Getzler proved (0.7) and (0.9) by essentially
the same
method.

\bigskip

\centerline{\bf \S 1. Moduli spaces}

\bigskip

In this section, we prove the Theorem 0.3.1 following the three step
procedure described in 0.6.

\medskip

{\bf 1.1. Marked trees and strata.} A tree is called {\it stable}
if $|v|\ne 2$ for all vertices $v$. If $|v|=1$ we call $v$
end vertex. Let $V_{\tau}^1$ be the set of end vertices.
An $n$--marking of $\tau$ is a bijection
$\mu :V^1_{\tau}\to\{ 1,\dots ,n\} .$ We also put
$V^0_{\tau}=V\setminus V^1_{\tau}$ and refer to it
as the set of interior vertices.

\smallskip

Let now $(C;x_1,\dots ,x_n)$ be a compact connected curve of
arithmetical genus zero with $n\ge 3$ labelled non--singular
points. The combinatorial structure of this curve is described
by the following stable tree with $n$--marking $(\tau ,\mu )$:
$V^0_{\tau}=\ \{$irreducible components of $C\}$,
$V^1_{\tau}=\{x_1,\dots x_n\} $; $\mu :x_i\mapsto i$;
an edge connects two interior vertices if the respective
components of $C$ have non--empty intersection; an edge
connects an interior vertex to an end vertex if the respective
point belongs to the respective component.

\smallskip

Denote now by $M(\tau ,\mu )\subset \overline{M}_{0,n}$
the set of points parametrising stable curves of the type
$(\tau ,\mu ).$ If $\tau$ has only one interior vertex,
$M(\tau ,\mu ):=M_{0,n}$ is the big cell. The following
statement summarises the main properties of these sets;
for a proof, see [Ke].

\smallskip

\proclaim{\quad 1.1.1. Proposition} a). $M(\tau ,\mu )$ is
a locally closed subset of $\overline{M}_{0,n}$ depending only
on (the isomorphism class of) $(\tau ,\mu ).$

b). $\overline{M}_{0,n}$ is the union of pairwise disjoint strata
$M(\tau ,\mu )$ for all marked stable $n$--trees $(\tau ,\mu ).$

c). For any $(\tau ,\mu )$,
$$
M(\tau ,\mu )\cong \prod_{v\in V^0_{\tau}}M_{0,|v|}.
$$
\endproclaim

\smallskip

Notice that there exists exactly one stable tree $\bullet$------$\bullet$
which does not correspond to any stable curve.

\smallskip

We can now calculate Poincar\'e polynomials.

\smallskip

\proclaim{\quad 1.1.2. Proposition} We have
$$
P_{M(\tau ,\mu )}(q)=\prod_{v\in V^0_{\tau}}P_{M_{0,|v|}}(q),
\eqno{(1.1)}
$$
$$
P_{M_{0,k}}(q)={{q^2-2}\choose{k-3}}(k-3)!. \eqno{(1.2)}
$$
\endproclaim

\smallskip

{\bf Proof.} (1.1) follows from the Proposition 1.1.1 and the
multiplicativity of Poincar\'e polynomials.

\smallskip

To prove (1.2), one can use the following geometric facts. First,
the morphism $\pi :\ \overline{M}_{0,n+1}\to \overline{M}_{0,n}$
forgetting the last marked point is (canonically isomorphic to)
the universal curve. Second, the infinity of the source consists of structure
sections and fibers at infinity of the target. Therefore, over the big cell
$M_{0,n}$ this morphism is a Zariski locally trivial fibration
with fiber $\bold{P}^1$, and $M_{0,n+1}=\pi^{-1}(M_{0,n})\setminus
\{${\it union of structure sections}$\}$.

\smallskip

{}From the addivity of Poincar\'e polynomials it follows that
$$
P_{M_{0,n+1}}(q)=P_{M_{0,n}}(q)P_{\bold{P}^1}(q)
-nP_{M_{0,n}}(q)=(q^2+1-n)P_{M_{0,n}}(q).
$$
Since $P_{M_{0,3}}(q)=1$, we get (1.2).

\medskip

Summarizing, we have for $n\ge 3$:
$$
P_{\overline{M}_{0,n}}(q)t^n=\sum_{{(\tau ,\mu )/(iso)}\atop{|V^1_{\tau}|=n}}
\prod_{v\in V^0_{\tau}}{{q^2-2}\choose{|v|-3}}(|v|-3)!
\prod_{v\in V^1_{\tau}}t,
\eqno{(1.3)}
$$
where $t$ is a new formal variable, and the sum is taken over
$n$-marked stable trees.

\medskip

{\bf 1.2. Passage to the standard marking.} Comparing (1.3) to
(0.22) and (0.23) we are more or less compelled to choose
$A=\{ *\}$ (one element set), $g^{**}=1,\ C_*=t,\ C_{**}=0$
(this gives weight zero to non--stable trees), and finally,
denoting by $C_k$ the component with $k\ge 3$ indices,
$$
C_k={{q^2-2}\choose{k-3}}(k-3)! .\eqno{(1.4)}
$$
In particular, we can forget about $f:\ F_{\tau}\to \{ *\}$.

\smallskip

This makes the weight of $(\tau ,\mu )$ depend only on $\tau /(iso)$,
but not $\mu$. Now, if $|V^1_{\tau}|=n$, the set of all $n$--markings of
$\tau$ consists of $n!$ elements and is effectively acted upon
by the group $\roman{Aut}\ \tau$. Therefore,
$$
\roman{card}\ \{ (\tau ,\mu )\}/(iso) =
\frac{n!}{|\roman{Aut}\ \tau|}\roman{card}\ \{ \tau\} /(iso).
$$
Putting together (1.3), (0.22), and (0.23), we see finally
that $\Phi (q,t)=Z^t$ where
$$
\Phi (q,t):= \frac{t^2}{2!}+\sum_{n\ge 3}\frac{t^n}{n!}
P_{\overline{M}_{0,n}}(q), \eqno{(1.5)}
$$
$$
Z^t:=\sum_{\tau /(iso)}\frac{1}{|\roman{Aut}\ \tau |}
\prod_{v\in V_{\tau}}C_{|v|}.
\eqno{(1.6)}
$$
The summation in (1.6) is now taken over all trees, the term
$t^2/2$ in (1.5) comes from the two--vertex tree, and the
generating function argument $t$ in (1.5) corresponds
precisely to the deformation parameter $t$ introduced at the end
of the subsection 0.6.

\smallskip

We will now use (0.26) in order to calculate
$$
\frac{\partial Z^t}{\partial t}=\frac{\partial \Phi (q,t)}{\partial t}:=
\varphi (q,t).
$$

\smallskip

{\bf 1.3. Potential.} From (0.24) and (1.4) one sees that
$$
S^t(\varphi )=-\frac{\varphi^2}{2}+t\varphi +\sum_{k\ge 3}C_k\varphi^k=
$$
$$
-\frac{\varphi^2}{2}+t\varphi +\sum_{k\ge 3}{{q^2-3}\choose{k-3}}
\frac{\varphi^k}{k(k-1)(k-2)}.
$$
This can easily be summed. We need only the derivative.

\medskip

\proclaim{\quad 1.3.1. Proposition} For generic $q$ we have
$$
\frac{\partial}{\partial \varphi}S^t(\varphi )=
\frac{(1+\varphi )^{q^2}-1-q^4\varphi}{q^2(q^2-1)}+t,
\eqno{(1.7)}
$$
and for $q=-1$,
$$
\frac{\partial}{\partial \varphi}S^t(\varphi )=
(1+\varphi )\ \roman{log}(1+\varphi )-2\varphi +t.
\eqno{(1.8)}
$$
\endproclaim

\medskip

{\bf 1.4. End of the proof.} We see now that (0.7), resp. (0.9),
are equations for the critical point $d_{\varphi}S^t=0$.
Differentiating them in $t$ and eliminating $(1+\varphi )^{q^2},$
resp $\roman{log}\ (1+\varphi ),$ we get (0.8), resp. (0.10).

\medskip

{\bf 1.5. Ramification of $\varphi (q,t)$ as a function of $t$.}
If $q^2$ is rational but $\ne 1,$ we see from (0.7) that $\varphi$
is an algebraic function of $t$ of genus 0. Otherwise it is
transcendental and infinitely valued. In order to understand
its topology, we can use the following classical trick.

\smallskip

Consider the differential equation for a function $y=y(x)$:
$$
yy_x=ax+by;\ a,b\in \bold{C}. \eqno{(1.9)}
$$
Let $w_{1,2}$ be roots of its {\it characteristic equation}
$$
w^2-bw-a=0. \eqno{(1.10)}
$$
Assume that $w_1\ne w_2$ and put
$$
A_1=\frac{w_1}{w_2-w_1},\ A_2=\frac{w_2}{w_1-w_2}
\eqno{(1.11)}
$$
so that $A_1+A_2=-1.$ A direct calculation shows:

\smallskip

\proclaim{\quad Proposition 1.5.1.} Put $w(x)=y(x)/x.$ Then the
general solution of (1.9) is given by the implicit equation
$$
Cx=(w-w_1)^{A_1}(w-w_2)^{A_2}, \eqno{(1.12)}
$$
where $C$ is an arbitrary constant.
\endproclaim

\smallskip

We can apply this to (0.8) putting
$$
y=1+q^2t-q^2\varphi,\ x=q^2t+q^2+1.
$$
Then we find
$$
w_1=1,\ w_2=q^{-2}, A_1=\frac{q^2}{1-q^2},\ A_2=\frac{1}{q^2-1}.
$$
One can calculate $C$ evaluating (1.12) at the point $t=0$
where we have $x=q^2+1,\ y=1,\ w=(q^2+1)^{-1}.$

\bigskip

\centerline{\bf \S 2. Configuration spaces}

\bigskip

In this section, we prove Theorems 0.4.1 and 0.4.2.

\medskip

{\bf 2.1. Nests and strata.} Let $X$ be a smooth compact algebraic variety.
The configuration space $X[n],\ n\ge 2,$ is defined in [FMPh]
as the closure of its big cell $X^n\setminus(\cup_{i<j}\Delta_{ij})$\quad
($\Delta_{ij}$ is the diagonal $x_i=x_j$) in
$X^n\times\prod_S\widetilde{X}^S,$ where $S$ runs over
subsets $S\subset\{ 1,\dots ,n\},\ |S|\ge 2;\ X^S$ denotes
the respective partial product of $X'$s, and $\widetilde{X}^S$
is the blow up of the small diagonal $\Delta_S$ in $X^S$.

\smallskip

Every $S$ determines a divisor at infinity $D(S)\subset X[n].$
Namely, let $\pi_S:\ X[n]\to X^S$ be the canonical projection.
Then $\pi_S^{-1}(\Delta_S)=\cup_{T\supset S}D(T).$

\smallskip

The natural stratification of $X[n]$ described in [FMPh]
consists of (open subsets of) intersections
$\overline{X(\Cal{S})}=\cap_{i=1}^{r}D(S_i)$ corresponding to sets
$\Cal{S}=\{ S_1,\dots ,S_r\}$ of subsets in $\{ 1,\dots ,n\}$
called {\it nests.}

\medskip

\proclaim{\quad 2.1.1. Definition} a). $\Cal{S}=\{ S_1,\dots ,S_r\}$
is a nest (or $n$--nest) if $|S_i|\ge 2$ for all $i$, and either
$S_i\subset S_j$ or $S_j\subset S_i$ for all $i,j$ such that
$S_i\cap S_j\ne \emptyset.$

\smallskip

In particular, $\Cal{S}=\emptyset$ is a nest, and
$\Cal{S}=\{ S\}$ is a nest, if $|S|\ge 2.$

\smallskip

b). A nest $\Cal{S}$ is called whole (resp. broken) if
$\{1,\dots ,n\}\in\Cal{S}$ (resp. $\{ 1,\dots ,n\}\notin\Cal{S}$).
\endproclaim

\smallskip

Denote by $X(\Cal{S})\subset\overline{X(\Cal{S})}=
\cap_{S\in\Cal{S}}D(S)$ the subset of points not belonging to smaller
closed strata. The following facts are proved in [FMPh].

\medskip

\proclaim{\quad 2.1.2. Proposition} a). For any $n\ge 2$ and
$n$--nest $\Cal{S},\ X(\Cal{S})$ is a locally closed subset of $X[n].$

\smallskip

b). $X[n]$ is the union of pairwise disjoint strata $X(\Cal{S})$
for all $n$--nests $\Cal{S}.$
\endproclaim

\medskip

{\bf 2.2. From nests to marked trees.} As in 1.1 we consider
a bijection $\mu :\ V^1_{\tau}\to\{ 1,\dots ,n\}$ as a part
of the appropriate marking for our problem. The
remaining data is supplied by choosing {\it orientation of
all edges.}

\smallskip

\proclaim{\quad 2.2.1. Definition} A tree $\tau$ marked in this
way is called admissible iff:

\smallskip

a). Every vertex of $\tau$ except of one has exactly one incoming
edge.

\smallskip

b). The exceptional vertex has only outgoing edges, and their number
is $\ge 2.$ This vertex is called source.

\smallskip

c). All interior vertices with possible exception of source
have valency $\ge 3.$
\endproclaim

\smallskip

\proclaim{\quad 2.2.2. Proposition} The following maps are (1,1):
$$
\{ broken\ n-nests\}\to\{ whole\ n-nests\}\to
\{admissible\ marked\ n-trees\} /(iso),
$$
$$
\Cal{S}\mapsto \Cal{S}\cup\{\{1,\dots ,n\}\}\mapsto
\tau{(\Cal{S})}=\tau{(\Cal{S}\cup\{\{ 1,\dots ,n\}\} )}.
$$

Here $\tau$ is defined by its sets of vertices and edges:
if $\Cal{S}=\{ S_1,\dots ,S_r\}$, then
$$
V_{\tau}=\{\widetilde{S}_1,\dots ,\widetilde{S}_{n+r}\}:=
\{S_1,\dots ,S_r,\{ 1\},\dots ,\{ n\}\},
$$
and an edge oriented from $\widetilde{S}_i$ to $\widetilde{S}_j$
connects these two vertices iff $\widetilde{S}_j\subset
\widetilde{S}_i$ and no $\widetilde{S}_k$ lies strictly
in between these two subsets.
\endproclaim

\smallskip

This is proved by direct observation. The following facts
are worth mentioning.

\smallskip

a). $\{ 1,\dots ,n\}$ is the source of $\tau (\Cal{S})$ for
any $\Cal{S}$.

\smallskip

b). $\{ 1\},\dots ,\{ n\}$ are all end vertices.

\smallskip

c). $i\in S_j$ iff one can pass from $S_j\in V_{\tau}$
to $\{ i\}\in V_{\tau}$ in $\tau$ by going always in positive
direction.

\smallskip

A reader is advised to convince him-- or herself that the source
has valency $\ge 2$ and all other interior vertices have
valency $\ge 3.$

\smallskip

Denote the source by $s$ and the set of the remaining interior
vertices $V^0_{\tau}.$

\medskip

\proclaim{\quad 2.2.3. Proposition ([FMPh])} The virtual
Poincar\'e polynomials of strata $X(\Cal{S})$ are given
by the following formulas (we add a formal variable $t$).

\smallskip

If $\Cal{S}$ is a broken $n$--nest, $s\in V_{\tau (\Cal{S})}$:
$$
t^nP_{X(\Cal{S})}(q)=
{{P_X(q)}\choose{|s|}}|s|!\times
\prod_{v\in V^0_{\tau (\Cal{S})}}
\kappa_m{{q^{2m}-2}\choose{|v|-3}}(|v|-3)!\times
\prod_{v\in V^1_{\tau (\Cal{S})}}t.
\eqno{(2.1)}
$$

If $\Cal{S}$ is a whole $n$-nest:
$$
t^nP_{X(\Cal{S})}(q)=P_X(q)\kappa_m{{q^{2m}-2}\choose{|s|-2}}
(|s|-2)!\times
$$
$$
\prod_{v\in V^0_{\tau (\Cal{S})}}
\kappa_m{{q^{2m}-2}\choose{|v|-3}}{(|v|-3)!}\times
\prod_{v\in V^1_{\tau (\Cal{S})}}t.
\eqno{(2.2)}
$$
\endproclaim

\smallskip

Comparing (2.1) and (2.2) one sees that one can express the
joint contribution of two nests corresponding to an admissible
marked tree $\tau$ as a product of local weights corresponding
to all vertices of $\tau$. The local weight of the source will be
$$
{{P_X(q)}\choose{|s|}}|s|!+P_X(q)\kappa_m
{{q^{2m}-2}\choose{|s|-2}}(|s|-2)!
$$
and the remaining local weights in (2.1) and (2.2) coincide and depend only
on the valency.

\medskip

{\bf 2.3. Passage to the standard marking.} We make the following
choices.

\smallskip

Put $A=\{ +,-\}.$ Interpret a mark $+$ (resp. $-$) on a flag as
incoming (resp. outgoing) orientation of this flag. Thus,
$f:\ F_{\tau}\to A$ is a choice of orientation of all flags.

\smallskip

Put $g^{+-}=g^{-+}=1,\ g^{++}=g^{--}=0.$ This makes the standard
weight of $(\tau ,f)$ vanish unless all {\it edges} are
unambiguously oriented by $f.$

\smallskip

Put $C_{+}=t$ (see (2.1) and (2.2)) and $C_{-}=0.$ The last
choice makes the standard weight vanish unless all end edges are
oriented outwards.

\smallskip

Put $C_{+-}=C{-+}=0.$ This excludes vertices of the type
$\to\bullet\to .$

\smallskip

Put also $C_{a_1,\dots ,a_k}=0$ if $\{ +,+\}\subset
\{ a_1,\dots ,a_k\} .$ This eliminates vertices with $\ge 2$
incoming edges.

\smallskip

For tensors with $k\ge 2$ minuses among the indices we put
$$
C_{-\dots -}={{P_X(q)}\choose{k}}k!+\kappa_mP_X(q)
{{q^{2m}-2}\choose{k-2}}(k-2)!
\eqno{(2.3)}
$$
(because only the source has all outgoing edges), and
$$
C_{+-\dots -}=\kappa_m{{q^{2m}-2}\choose{k-2}}(k-2)!
\eqno{(2.4)}
$$
(cf. (2.1) and (2.2)).

\smallskip

The standard weight of a marked tree defined by this data again is
independent on the part $\mu :\ V^1_{\tau}\to\{ 1,\dots ,n\}$
of the initial marking which accounts for the factor
$\dfrac{n!}{|\roman{Aut}\ \tau |}$ below.

\smallskip

Summarizing, we put
$$
\Phi_X(q,t):=\sum_{n\ge 2}\frac{t^n}{n!}P_{X[n]}(q),
\eqno{(2.5)}
$$
$$
Z^t:=\sum_{\tau /(iso)}\frac{1}{|\roman{Aut}\ \tau |}
\sum_{f:F_{\tau}\to\{ +,-\}}\prod_{\alpha\in E_{\tau}}
g^{f(\partial \alpha )}\prod_{v\in V_{\tau}}C_{f(\sigma v)},
\eqno{(2.6)}
$$
and get from the previous discussion
$$
Z^t=\Phi_X(q,t),\qquad \frac{\partial}{\partial t}Z^t:=
\phi_X(q,t).  \eqno{(2.7)}
$$

\medskip

{\bf 2.4. Potential.} We change notation: $\varphi_+=x\ ,
\varphi_-=y.$ From 2.3 we see that (already $t$-deformed)
potential is
$$
S^t(x,y)=-xy+tx+\kappa_m\sum_{k=2}^{\infty}
{{q^{2m}-2}\choose{k-2}}\frac{xy^k}{k(k-1)}+
$$
$$
\sum_{k=2}^{\infty}{{P_X(q)}\choose{k}}y^k+\kappa_mP_X(q)
\sum_{k=2}^{\infty}{{q^{2m}-2}\choose{k-2}}{\frac{y^k}{k(k-1)}}
\eqno{(2.8)}
$$
(we have two arguments $x,y$ but only one $t=t_+$ because
$C_-=0$).

\smallskip

We must solve the system
$$
\frac{\partial S^t}{\partial x}|_{x^0,y^0}=
\frac{\partial S^t}{\partial y}|_{x^0,y^0}=0,  \eqno{(2.9)}
$$
and (0.26) then tells us that
$$
\frac{\partial}{\partial t}Z^t=\varphi_X(q,t)=x^0 . \eqno{(2.10)}
$$

Again, $S^t(x,y)$ can be easily summed. To write down the functional
equation, we need only $x$--derivative which for general $q$ is
$$
\frac{\partial S^t}{\partial x}=-y+t+
\kappa_m\frac{(1+y)^{q^{2m}}-1-q^{2m}y}{q^{2m}(q^{2m}-1)}.
\eqno{(2.11)}
$$
For $q=-1$:
$$
\frac{\partial S^t}{\partial x}=-y+t+m[(1+y)\ \roman{log}(1+y)-y].
\eqno{(2.12)}
$$
\medskip

{\bf 2.5. End of the proof.} We now see that (0.15), resp (0.18),
are the equations defining $y^0$. Differentiating
in $t$ we get (0.16) and (0.19). And since $S^t(x,y)$ is linear
in $x$, the vanishing of the $y$--derivative gives an explicit
expression of $x^0$ via $y^0$:
$$
\varphi_X(q,t)=P_X(q)
\frac{(1+y^0)^{P_X(q)}+(q^{2m}+\kappa_m-1)y^0-q^{2m}t-1}{1+(1-q^{2m}-
\kappa_m)y^0+q^{2m}t}.
$$
To see that this is equivalent to (0.17) one can differentiate (0.17)
in $t$ and use (0.16).

\medskip

{\bf 2.6. Ramification of $y^0$.} Replaying the game of 1.5,
we put (changing the meaning of $x,y$ in favor of those in 1.5):
$$
y=y(q,t):=q^{2m}t+1-(q^{2m}+\kappa_m-1)y^0(q,t),
$$
$$x:=t+\frac{q^{2m}+\kappa_m}{q^{2m}},\qquad w(q,t)=y/x.$$
Then (0.16) becomes
$$
yy_x=-q^{2m}x+(q^{2m}+1)y
$$
so that in the notation of 1.5
$$
w_1=1,\ w_2=q^{2m},\ A_1=\frac{1}{q^{2m}-1},\ A_2=
\frac{q^{2m}}{1-q^{2m}},
$$
and finally
$$
Cx=(w-w_1)^{A_1}(w-w_2)^{A_2}
$$
for some $C.$

\bigskip

\centerline{\bf \S 3. Multiple coverings}

\bigskip

{\bf 3.1. Kontsevich's formula for Problem C.} Kontsevich represents
$m_d$ as a rational function of two variables $\lambda_1,\lambda_2$
which is formally homogeneous of degree zero and actually is expected to be
a constant.

\smallskip

Geometrically, this statement must be a corollary of Bott's fixed
point formula for {\it smooth stacks}. The $\lambda$--variables
in this context are coordinates of a toric vector field on the
target $\bold{P}^1.$ Until this has been worked out, we simply
go ahead with Kontsevich and take this independence for granted.

\smallskip

The function in question is a sum of contributions indexed
by isomorphism classes of connected trees $\tau$ endowed with
markings: each vertex $v$ is marked by $f_v=1$ or 2 so that no neighbors
have the same mark; each edge $\alpha$ is marked by a positive
integer $d_{\alpha}.$ Only those marked trees contribute to $m_d$
for which $\roman{deg}\ \tau :=\sum_{\alpha}d_{\alpha}=d.$

\smallskip

We introduce the following notation for a marked tree $\tau$:\
$F$=\ the number of vertices marked by 2; $\sigma_v=\sum_{\alpha :
v\in \alpha}d_\alpha ;$
$w_i=\sum_{v:f_v=i}(|v|-1),\ i=1,2.$

\smallskip

Then we have
$$
m_d=(\lambda_1-\lambda_2)^{2-2d}\sum_{\tau :\roman{deg} \tau =d}
\frac{1}{|\roman{Aut}\ \tau |}(-1)^{d+F}\lambda_1^{2w_1}\lambda_2^{2w_2}
V(\tau )E(\tau ),
$$
$$
V(\tau )=\prod_v \sigma_v^{|v|-3},\
E(\tau )=\prod_{\alpha}\frac{d_{\alpha}^3}{d_{\alpha}!^2}
\prod_{a+b=d_\alpha ;a,b\ge 1}(a\lambda_1+b\lambda_2)^2.
$$

\proclaim{\ 3.2. Theorem} $m_d=d^{-3}.$
\endproclaim

{\bf Proof.} We will calculate the value of $m_d$ at $\lambda_1=1,\
\lambda_2=0.$ The drastic simplification results from the fact
that the factor $\lambda_2^{2w_2}$ vanishes unless $w_2=0$.
Now, $w_2=0$ implies that $\tau$ has no vertices of multiplicity
$\ge 2$ marked by 2. Hence $\tau$ either has only one edge,
or is a star with central vertex marked 1, and end vertices marked
by 2. We will consider the first case as one ray star as well.

\smallskip

Now, let $\tau $ be such a star of degree $d.$ The set
$\{ d_{\alpha}\}$ forms a partition of $d$ into positive
summands which uniquely defines the isomorphism class of $\tau .$
It is convenient to write this partition as the set
of multiplicities $R=\{r_1,r_2,\dots\} ,$ where $r_i$= the number
of edges marked by $i$
so that $\sum_iir_i=d.$ Obviously, $|\roman{Aut}\ \tau |=\prod_ir_i!.$

\smallskip

After some reshuffling, our assertion thus reduces to the following
identity:
$$
(?)\qquad \sum_R \frac{1}{\prod_ir_i!}\prod_i(-\frac{d}{i})^{r_i}
=(-1)^d.
$$
Now, the left hand side of (?) can be obtained in the following way.
Consider the formal series
$e^{\sum_{i\ge 1}y_it^i}$, take its terms of degree $d$ in $t$
and put in them $y_i=-d/i.$ But we can clearly proceed in reverse
order first making the substitution $y_i=-d/i.$ Then
the series in the exponent becomes $\sum_i(-d/i)t^i=d \roman{log}(1-t),$
so that finally we get the coefficient of $t^d$ in $(1-t)^d.$ QED

\smallskip

{\it Remark.} One can observe that $(-1)^d$ coincides with the contribution
of just one trivial partition: $r_d=1$. The remaining terms cancel.
Geometrically, this means that degenerating configurations do
not contribute with this choice of vector field. Algebraically,
this can be rewritten as an equality of two sums, one over
proper partitions with odd, another with even number of
summands.

\bigskip

\centerline{\bf References}

\bigskip

[AM] P. S. Aspinwall, D. R. Morrison. {\it Topological field
theory and rational curves.} Comm. Math. Phys. 151 (1993), 245--262.

\smallskip

[BG] A. Beilinson, V. Ginzburg. {\it Infinitesimal structure
of moduli spaces of $G$-bundles.} Int. Math. Res. Notes,
4 (1992), 63--74.

\smallskip

[FMPh] W. Fulton, R. MacPherson. {\it A compactification
of configuration spaces.} Ann. of Math. 139 (1994), 183--225.

\smallskip

[GK] V. Ginzburg, M. Kapranov. {\it Koszul duality for operads.}
Preprint, 1993.

\smallskip

[Ke] S. Keel. {\it Intersection theory of moduli spaces
of stable $n$--pointed curves of genus zero.} Trans.
AMS, 330 (1992), 545--574.

\smallskip

[Ko] M. Kontsevich. {\it Enumeration of rational curves via
torus action.} Preprint MPI/94--39.

\enddocument